\documentclass[journal]{IEEEtran}

\usepackage[final]{graphicx}
\graphicspath{{figs}}
\usepackage{tikz}
\ifCLASSINFOpdf
\else
\fi
\usepackage[cmex10]{amsmath}
\usepackage{bm}
\usepackage{amssymb}
\usepackage{xcolor}
\usepackage{cite}
\usepackage{soul}
\usepackage[caption=false]{subfig}

\newcommand{\br}{\bar{r}}

\begin{document}
\title{Time-dependent Debye-Mie Series Solutions for Electromagnetic Scattering}
\author{
		  Jie Li,~\IEEEmembership{Student Member,~IEEE,}
		  Balasubramaniam~Shanker,~\IEEEmembership{Fellow,~IEEE}
		  \thanks{J. Li is with the Department of Electrical and Computer Engineering, Michigan State University, East Lansing, MI 48824-1226 email: jieli@egr.msu.edu}
		  \thanks{B. Shanker is with the Department of Electrical and Computer Engineering, Michigan State University}}
\maketitle

\begin{abstract}

Frequency domain Mie solutions to scattering from spheres have been used for a long time. However, deriving their transient analogue is a challenge as it involves an inverse Fourier transform of the spherical Hankel functions (and their derivatives) that are convolved with inverse Fourier transforms of spherical Bessel functions (and their derivatives). Series expansion of these convolutions are highly oscillatory (therefore, poorly convergent) and unstable. Indeed, the literature on numerical computation of this convolution is very sparse. In this paper, we present a novel quasi-analytical approach to computing transient Mie scattering that is both stable and rapidly convergent. The approach espoused here is to use vector tesseral harmonics as basis function for the currents in time domain integral equations together with a novel addition theorem for the Green's functions that renders these expansions stable. This procedure results in an orthogonal, spatially-meshfree and singularity-free system, giving us 
a set of one dimensional Volterra Integral equations. Time-dependent multipole coefficients for each mode are obtained via a time marching procedure.  Finally, several numerical examples are presented to show the accuracy and stability the proposed method.

\end{abstract}


\section{Introduction}

Debye-Mie series solution \cite{Mie1908,Stratton1941,Ishimaru1990} is one of the most useful tools in time-harmonic analysis of electromagnetic (EM) scattering, and has found in extensive use in both optics and electromagnetics \cite{Kerker1969,Hergert2012}. Their applications range from fields as diverse as light scattering from small particles \cite{Mishchenko2002} to combustion \cite{Guelder2000} to analysis of small antennas \cite{Schuller2009} to  weather \cite{Deirmendjian1964} to validation and calibration computational and measurement setups to photonics \cite{Rybin2009} to biological applications \cite{Bassan2010}. However, analytical solutions to scattering from a sphere only exist in the Fourier domain. Equivalent literature for transient analysis is very sparse. While it is indeed possible to to obtain transient scattering response via Fourier transforms,  the challenge of obtaining a stable direct transient solution has remained unsolved for several decades. While attempts to arrive at a 
solution 
have been made \cite{Buyukdura1997} via explicit Fourier transforms of Bessel and Hankel 
functions, these approaches are highly oscillatory and therefore, unstable. Aside from an intellectual challenge, the development of methods to solve for scattering from sphere can find applications in a number of areas; from scattering from collection of particles to study of non-linear response of gain particles \cite{Liu2003} to  imaging \cite{Jain2006} to applications to dielectric resonator antennas \cite{Leung1993}, and so on. 

While direct time-domain Mie series like solutions for scattering do not exist, there has been extensive interest in developing multipole theory for representations of radiated fields due to spatially bounded sources. Time domain multipole theory has found applications in both acoustics as well as electromagnetics \cite{Shlivinski2001,Heyman1996,Granzow1966,Davidon1973,Marengo1998,Devaney1974}, with applications that range from radiation due to volumetric sources \cite{Heyman1996a} to near field scanning \cite{Shlivinski1999} to near-to-far field transforms \cite{Oetting2005}. However, multipole representations require higher order derivatives in time. As a result, there is a paucity of numerical data and the treatment is largely theoretical. This comment does not apply to fast solvers that use a  
Whittaker 
representation \cite{Whittaker1903} of the radiated field, and there are substantial numerical results in both acoustics \cite{Ergin2000} and electromagentics \cite{Shanker2003}. 

However, compared with the volume of literature on representation of transient fields due to quiescent sources, research on transient scattering analysis using spherical harmonics is somewhat limited. The methods investigated thus far fall into two categories; (i) using inverse Fourier/Laplace transforms of Bessel and Hankel functions within a mode matching framework,  and (ii) using a novel representation of the retarded potential Green's function.
The first time-domain Mie algorithms were presented for the scalar Helmholtz equation in \cite{Buyukdura1997}. It is evident from the work that obtaining a stable algorithm is difficult as the inverse Fourier transform  of the spherical Hankel function is a Gagenbauer polynomial that grows as a function of time. The inverse Fourier transform in the work has to be understood in the convolution sense because the actual inverse Fourier transform doesn't exist. A more recent paper examined the causes of instability, errors in the extent literature, and presented a recursive convolution approach that is stable for both late time and high order harmonics \cite{Greengard2014} for scalar fields. Extension of classical Fourier methods to electromagnetic fields is significantly more complicated \cite{Azizoglu2004}, and given the instability for scalar fields, it is not clear that this will be stable for vector fields. It is possible that the method developed in \cite{Greengard2014} may be extended to the vector case.
The second method for solving scattering from spheres was developed before \cite{Greengard2014}, but took an entirely different approach \cite{Li2014}. The underlying thesis of this approach was to utilize the fact that if the trace of quantities of interest were to be represented in terms of tesseral harmonics, then the solution of an integral equation composed of these trace quantities would yield the Mie solution. This statement can be readily proven in the frequency domain \cite{Tai1994}, and the proofs arise from a spherical harmonic representation of the Green's 
function. In time domain, this would translate to a representation of the retarded potential Green's functions, which would bring us \textit{back} the original set of difficulties. In \cite{Li2014}, we proposed a novel representation of the Green's function that obviates these difficulties, and produces stable and convergent results.     

The work presented herein would follow in the vein of \cite{Li2014} and leverage time domain integral equation (TDIE) methods to develop a quasi-analytical solution time dependent EM scattering from a sphere. The approach presented is based on  vector spherical harmonics expansion for the time-dependent dyadic Green's function (to be more specific, its tangential trace) of both electric and magnetic field types, and integrated the results presented in \cite{Li2014}. The use of tesseral harmonics for representing the spatial variation of the trace quantities on the surface of the sphere permits (i) a mesh free solution and (ii) analytical evaluation of the integral required for the TDIE, leaving only a Volterra integral in time. The specific contributions of the paper are as follows: We 
\begin{itemize}
	\item present a stable spherical expansion for time domain dyadic Green's functions of both electric and magnetic types.
        \item present a spatially mesh-less and singularity-free integral equations for EM scattering from spherical objects
        \item develop a set of reduced 1-D Volterra integral kernels that are solved with a higher order Galerkin scheme
	\item present several results demonstrating the accuracy, convergence and stability of the proposed methods.  
\end{itemize}

The remainder of the paper is organized as follows. Section \ref{sec:probs} describes the problem to be solved  and gives the integral equation 
with time-domain dyadic Green's function.  Section \ref{sec:form} derives the vector spherical harmonic expansion for the
tangential trace of the time-dependent dyadic Green's function.
Section \ref{sec:time_mie} gives a brief approach to solve the reduced integral equation to obtain the time-dependent spherical multipoles. 
Section \ref{sec:results} presents a number of results that validate the accuracy, convergence and stability of the presented method.  Finally, we summarize our contribution and outline future research directions in Section \ref{sec:conc}.

\section{Problem Statement}\label{sec:probs}
Consider a perfectly electrically conducting spherical object that occupies a volume $D_1 \subset \mathbb{R}^3$ residing in a homogeneous background medium occupying   {$D_0=\mathbb{R}^3\backslash D_1$}, characterized by permittivity $\epsilon_0$ and permeability $\mu_0$. The boundary of the scatterer is denoted using $\Omega_1 = \partial D_1$, and is equipped with an outward pointing normal $\hat{n} (\br)$ and a radius of $r=a$.  On the interface $\Omega_1$, the tangential electric field vanishes.

The electric field integral equation (EFIE) and magnetic field integral equation (MFIE), their linear combination, can be used to solve the electric current on the surface $\Omega_1$ (usually with spatial meshing). The EFIE and MFIE are usually written as follows, respectively.
\begin{equation}
\label{eq:TDEFIE}
 \hat{n} \times \hat{n} \times  \mathbf{\mathcal{L}}(\mathbf{J}(\bar{r},t))  = \hat{n} \times \hat{n} \times {\bf E}^i (\bar{r},t) 
\end{equation}
\begin{equation}
\label{eq:TDMFIE}
\mathbf{J}(\bar{r},t) + \hat {n} \times  \mathbf{\mathcal{K}}(\mathbf{J}(\bar{r},t))  = \hat {n} \times \mathbf{H}^{i}(\bar{r},t) 
\end{equation}
where operators $\mathbf{\mathcal{L}}(\cdot) $ and $\mathbf{\mathcal{K}}(\cdot)$ are spatial-temporal (four-dimensional) operators associated with dyadic Green's function.
\begin{equation}
\mathbf{\mathcal{L}}(\mathbf{X}(\bar{r},t)) =\mu \partial_t \int_{S'} \tilde{ \bf G}_{e0}(\bar{r},\bar{r}') \otimes \mathbf{X}(\bar{r},t) dS'
\end{equation}
\begin{equation}
\mathbf{\mathcal{K}}(\mathbf{X}(\bar{r},t)) = - \int_{S'}  \tilde {\bf G}_{m0}(\bar{r},\bar{r}') \otimes \mathbf{X}(\bar{r},t) dS' 
\end{equation}
In the above definitions, $ \tilde{ \bf G}_{e0}$ and $ \tilde{ \bf G}_{m0}$, respectively, denote dyadic Green's functions of electric type and magnetic type in free space.The solutions to the above equations are typically effected using a discrete representation of the scatterer, representing the current using appropriate spatial and temporal basis on this discrete geometric representation and finally, creating a set of linear equations using testing procedures in both space and time. But if the scatterer is spherical/canonical, one can take advantage of the geometry. This is done easily in the frequency domain. In what follows, its transient analogue is presented. First, the time-dependent Dyadic Green's functions (to be more specific, their tangential  traces) are first expanded using vector spherical harmonics. And then the TD-EFIE and TD-MFIE are reduced into a simple Volterra integral equations 
without any spatial meshing and singularity in the integral kernels. It permits mode-by-mode 
solutions for coefficients of time-domain Debye-Mie series.

\section{Formulation}\label{sec:form}

\subsection{Natural Spatial basis: Vector spherical harmonics}

On a sphere, vector spherical harmonics (VSH) form a natural basis set to represent the spatial variation of the current \cite{Hill1954,Carrascal1991}. Their relations to spherical harmonics are given as
\begin{equation}
\mathbf{\Psi}_n^m({\hat r}) = \mathbf{\Psi}_n^m(\theta,\phi) =\dfrac{r}{ \sqrt{n(n+1)}} \nabla^{t} Y_n^m(\theta,\phi) 
\end{equation}
\begin{equation}
\mathbf{\Phi}_n^m({\hat r})  = \mathbf{\Phi}_n^m(\theta,\phi) =  \dfrac{1}{ \sqrt{n(n+1)}} \bar {r}\times \nabla^{t} Y_n^m(\theta,\phi)
\end{equation}
where 
\begin{equation}
Y_n^m(\theta,\phi) = \sqrt{  \dfrac {2n+1}{4\pi}  \dfrac {(n-m)!}{(n+m)!}  }  P_n^m(\cos \theta)e^{jm\phi}
\end{equation}
In the above equations, $P_n^m$ is the associated Legendre function of degree $n$ and order $m$ of the mode. VSHs comprise a complete and orthogonal tangential basis set. Therefore, the current on the sphere can be written as 
\begin{equation}
 {\bf J}(\bar{r},t) = \sum_{n,m} J_{nm}^1(t) \mathbf{\Psi}_n^m(\theta,\phi) + J_{nm}^2 (t) \mathbf{\Phi}_n^m(\theta,\phi)
\end{equation}
Given the coefficients of both modes, one can get the scattered electric and magnetic fields  using 
\begin{equation}
\mathbf{E}^s(\bar{r},t) = -\mathbf{\mathcal{L}}(\mathbf{J}(\bar{r}',t))
\end{equation}
\begin{equation}
\mathbf{H}^s(\bar{r},t) = -\mathbf{\mathcal{K}}(\mathbf{J}(\bar{r}',t))
\end{equation}
The basis here used is different from the basis used as in time-harmonic Debye-Mie series, where the field, rather than the current, is represented using the vector spherical wave functions (VSWF) with radial dependence. However, we use VSWF next to derive the expansions for the time-dependent dyadic Green's functions in time-domain integral equations.

\subsection{Spherical expansion of Time-domain Dyadic Green's function}\label{sec:sph_exp}
The time-domain dyadic Green's function of electric type is the solution to the following wave equation
\begin{equation}
\nabla \times \nabla \times \tilde {\bf G}_{e0} (\bar{r},\bar{r}',t) - \dfrac{\partial_t^2}{c^2} \tilde {\bf G}_{e0}(\bar{r},\bar{r}',t) 
=  {\it \tilde I} \delta (\bar{r}-\bar{r}')\delta(t)
\end{equation}
and its dyadic form is notionally written as 
\begin{equation}
 \tilde {\bf G}_{e0} (\bar{r},\bar{r}',t) = (\tilde I - c^2 \partial_t^{-2} \nabla \nabla)\dfrac{\delta(t-\frac{r-r'}{c})}{4\pi R}
\end{equation}
where $\partial_t^{-2}$ stands for a two-fold integral with respect to time, ${\tilde I}$ is the idempotent, and $\dfrac{\delta(t-\frac{R}{c})}{4\pi R}$ is retarded potential. As shown later,  $\partial_t^{-2}$ is not evaluated explicitly. To glean more insight into the analysis, we start with the frequency domain counterpart of the dyadic Green's function:
\begin{equation}
\tilde {\bf G}_{e0}(\bar{r},\bar{r}',\omega)
=\Big[ {\it \tilde I} + \dfrac{1}{k^2} \nabla \nabla \Big] G_0(\bar{r},\bar{r}',\omega)
\end{equation}
where $k$ and $\omega$ denote the wavenumber and the angular frequency, respectively. The representation of the dyadic Green's function in spherical coordinates is well known, has been used extensively \cite{Tai1994}, and reads as  
\begin{equation}
\label{eq:dyadGF_exp}
\begin{split}
 \tilde {\bf G}_{e0}(\bar{r},\bar{r}',\omega)
= jk \sum_{n,m} \Big[  
  & {\bf N}_{nm}^{(4)}(\bar{r},k) {\bf N}_{nm}^{(1)\ast}(\bar{r}',k) \\
&  ~~~ +   {\bf M}_{nm}^{(4)}(\bar{r},k) {\bf M}_{nm}^{(1)\ast}(\bar{r}',k)
\Big]
\end{split}
\end{equation}
for $|\bar{r}| > |\bar{r}'|$. In the above equations, ${\bf N}_{nm}^{(p)}(k,\bar{r})$ and ${\bf M}_{nm}^{(p)}(k,\bar{r})$ are the TM$^r$ and TE$^r$ vector spherical wave functions (VSWF) of degree $n$ and order $m$, respectively. In what follows, we use $r = |\bar{r}|$ and $r' = |\bar{r}'|$. The explicit expressions for ${\bf N}_{nm}^{(m)}(k,\bar{r})$ and ${\bf M}_{nm}^{(p)}(k,\bar{r})$  can be found in the appendix in terms of spherical Bessel $z_n^{(1)}(kr')$ or spherical Hankel functions $z_n^{(4)}(kr)$.  As one can see from \eqref{eq:dyadGF_exp}, the expansion is in terms of angular dependent wave function (angular dependence is integrated with radial dependence at this moment), rather than the vector spherical harmonics. Furthermore, the frequency dependence is sprinkled throughout, and as a result taking an inverse Fourier transform is difficult.  In order to fully exploit the orthogonality of VSHs as well as make the expressions more amenable to taking an inverse Fourier 
transform, it is 
necessary to recast the Green's function in terms of its tangential trace. This can be effected via the following: Taking the tangential trace of the dyadic Green's function via
\begin{equation}
\begin{split}
\tilde {\bf G}_{e0}^{tt}(\bar{r},\bar{r}',\omega)
= \hat {r}  \times \hat {r} \times  \tilde {\bf G}_{e0} (\bar{r},\bar{r}',\omega)  \times \hat{r}' \times \hat{r}' 
\end{split}
\end{equation}
yields
\begin{equation}
\label{eq:dyadGF_trace_def}
\begin{split}
\tilde {\bf G}_{e0}^{tt}(\bar{r},\bar{r}',\omega)
= - jk   \hat {r}  \times \hat {r} \times \sum_{n,m} \Big[  
 {\bf N}_{nm}^{(4)}(k,\bar{r}) {\bf N}_{nm}^{(1)\ast}(k,\bar{r}') \\
+   {\bf M}_{nm}^{(4)}(k,\bar{r}) {\bf M}_{nm}^{(1)\ast}(k,\bar{r}')
\Big] \times \hat{r}' \times \hat{r}' \\
 =-jk \sum_{n,m} \Big[  
\dfrac{\big[ kr{z_n^{(4)}}(kr) \big]'}{kr}   \dfrac{\big[kr'{z_n^{(1)\ast}}(kr') \big]'}{kr'}   {\bf \Psi}_{nm} (\hat{r})   {\bf \Psi}_{nm} (\hat{r}') \\
+ z_n^{(4)}(kr) z_n^{(1)\ast}(kr') {\bf \Phi}_{nm}(\hat{r}){\bf \Phi}_{nm}(\hat{r}')
\Big]
\end{split}
\end{equation}
It can be shown that the trace of the Green's function has the following properties when operating on the tangential field ${\bf X}$ that is defined on a spherical surface:
\begin{subequations}
 \begin{equation}
\tilde {\bf G}_{e0}^{tt}(\bar{r},\bar{r}',\omega) \cdot{\bf X}
=   - \tilde {\bf G}_{e0} (\bar{r},\bar{r}',\omega) \cdot  {\bf X} 
\end{equation}
 \begin{equation}
{\bf X} \cdot   \tilde {\bf G}_{e0}^{tt}(\bar{r},\bar{r}',\omega) 
= - {\bf X} \cdot    \tilde {\bf G}_{e0} (\bar{r},\bar{r}',\omega) 
\end{equation}
\end{subequations}
As a result of these properties, one can use $\tilde {\bf G}_{e0}^{tt}(\bar{r},\bar{r}',\omega)$ instead of $\tilde {\bf G}_{e0} (\bar{r},\bar{r}',\omega)$. This property is valid in both frequency and time domain. As is evident from \eqref{eq:dyadGF_trace_def}, $\tilde {\bf G}_{e0}^{tt}(\bar{r},\bar{r}',\omega)$ is expressed purely in terms of VSHs, and consequently one can exploit the orthogonality of these functions as well as  take its inverse Fourier transform. The resulting expressions for the dyadic Green's function can be written as 
\begin{equation}
\label{eq:invFT_GF}
\begin{split}
& \tilde {\bf G}_{e0}^{tt}(\bar{r},\bar{r}',t)
= \\
 \sum_{n,m} \Big[ 
& {\cal F}^{-1}   \big[  jk \dfrac{\big[ kr{z_n^{(4)}}(kr) \big]'}{kr}  \dfrac{\big[kr'{z_n^{(1)\ast}}(kr') \big]'}{kr'}  \big]
  {\bf \Psi}_{nm} (\hat{r})   {\bf \Psi}_{nm} (\hat{r}')  \\
& ~~~~ + {\cal F}^{-1} \big[ jkz_n^{(4)}(kr) z_n^{(1)\ast}(kr')  \big] {\bf \Phi}_{nm}(\hat{r}){\bf \Phi}_{nm}(\hat{r}')
 \Big]\\
& = \sum_{n,m} \Big[ 
I_n^1  {\bf \Psi}_{nm} (\hat{r})   {\bf \Psi}_{nm} (\hat{r}') +I_n^2 {\bf \Phi}_{nm}(\hat{r}){\bf \Phi}_{nm}(\hat{r}')
 \Big]  
\end{split}
\end{equation}

It can be shown that $I_n^1$ in \eqref{eq:invFT_GF} can be rewritten as
\begin{equation}
\label{eq:invFT_GF_1stterm}
\begin{split}
& I_n^1 = \dfrac{ c^2\partial_t^{-2}}{rr'}\partial_r\partial_{r'} \Big[ r r' {\cal F}^{-1}   \big[ jk  {z_n^{(4)}}(kr) {z_n^{(1)\ast}}(kr') \big]  \Big] \\
& ~~~~~  = \dfrac{ c^2\partial_t^{-2}}{rr'}\partial_r\partial_{r'} \Big[ r r'   I_n^2  \Big] 
\end{split}
\end{equation}
As is evident from \eqref{eq:invFT_GF} and \eqref{eq:invFT_GF_1stterm}, one needs to obtain the explicit expression for {\em only} $I_n^2$
in order to evaluate an expression of transient electric dyadic Green's function. As is well known, the approach of using the inverse Fourier transforms of spherical Bessel and Hankel functions and then evaluating their convolution is unstable, because the inverse Fourier transform of spherical Hankel function doesn't exist. In the next section, we present an overview of a method that was recently developed to address this very problem \cite{Li2014}, and then use this to recover expressions for both the trace of the electric and magnetic dyadic Green's functions. 


The starting point for deriving an expression that approximates $I_2$ or provides an alternative representation arises from examining the representation for the Green's function for the Helmholtz equation, viz., 
\begin{equation}
\label{eq:addThrmHelmholtz}
\begin{split}
\dfrac{\exp{[-jkR}]}{4 \pi R} & = -jk\sum_{n,m} z^{(4)}_{n}(k r) z_{n}^{(1)}(k r') Y_n^m({\hat r}) {Y_n^{m}}^{*} ({\hat r}')\\
\dfrac{\delta(t - R/c)}{4\pi R} & = \sum_{n,m} {\cal F}^{-1} [ -j kz^{(4)}_{n}(k r) z_{n}^{(1)}(k r')  ]Y_n^m({\hat r}) {Y_n^{m}}^{*} ({\hat r}')
\end{split}
\end{equation}
While this is the well known classical representation, an alternate representation may be derived using 
\begin{equation}
f (x) = \sum_n a_n P_n (x);~~a_n = \dfrac{2n+ 1}{2}\int_{-1}^{1} dx~f(x) P_n (x)
\end{equation}
Using $x =\cos \gamma = (r^2 + {r'}^2 - R^2)/(2 r r')$ where $\gamma$ denotes the angle between the points ${\bf r}$ and ${\bf r}'$ or equivalently $R= \sqrt{r^2 + r'^2 - 2 rr' x}$, and $f (x) = \exp{[- j k R]} /(4 \pi R)$, it is trivial to show that 
\begin{equation}
\dfrac{\exp{[-jkR]} }{4 \pi R} 
= \sum_{n=0}^{\infty} \dfrac{2n+1}{8\pi rr'}
\int_{|r - r'|}^{r + r'} \exp [-j k R] P_n \left ( x \right )dR 
\end{equation}
The inverse Fourier transform to the above expression yields 
\begin{equation}
\label{eq:altaddThrm}
\begin{split}
& \dfrac{\delta(t - R/c)}{4\pi R}  = \sum_{n=0}^{\infty} \dfrac{(2n + 1)c}{8\pi  r r'} 
           P_n\left (\frac {\xi} {2 r r'} \right ) \mathcal{P}_{\alpha \beta}(t) P_n( x ) \\
           & = \sum_{n,m} \dfrac{c}{2\pi  r r'} 
           P_n\left (\frac {\xi} {2 r r'} \right ) \mathcal{P}_{\alpha \beta}(t) Y_n^m({\hat r}) Y_n^{*m} ({\hat r}')
\end{split}
\end{equation}
where $\xi = r^2 + r'^2 -c^2 t^2$, $\mathcal{P}_{\alpha \beta}(t) $ is a pulse function (the value is $1$ when $t \in [\alpha,\beta]$, and zero otherwise) with $\alpha = |r - r'|/c $ and $ \beta =(r + r')/ {c}$. In obtaining  these expressions, we have used the addition theorem for Legendre polynomials. Comparing \eqref{eq:altaddThrm} to \eqref{eq:addThrmHelmholtz}, it is apparent that $I_n^2$ takes the form
\begin{equation}
\begin{split}
\label{eq:FT_I2}
& {\cal F}^{-1} \big[ - jkz_n^{(4)}(kr) z_n^{(1)\ast}(kr')  \big]  \\
& = \dfrac{c}{2rr'}P_n\Big(\dfrac{\xi}{2rr'}\Big)P_{\alpha \beta}(t) 
\equiv K_n^{(0)}(r,r',t) 
\end{split}
\end{equation} 
Using \eqref{eq:FT_I2} together with \eqref{eq:invFT_GF_1stterm}, one can obtain the dyadic Green's function of the electric type as follows: 
\begin{equation}
\begin{split}
& \tilde {\bf G}^{tt}_{e0}(\bar{r},\bar{r}',t) 
 =  \sum_{n,m} \Big[ K_n^{(0)}(r,r',t)  {\bf \Phi}_{nm}(\hat{r}){\bf \Phi}_{nm}(\hat{r}') \\
& ~~ + \dfrac{c^2\partial_t^{-2}}{rr'}\partial_r\partial_{r'} \Big [  rr'K_n^{(0)}(r,r',t) \Big] {\bf \Psi}_{nm} (\hat{r})   {\bf \Psi}_{nm} (\hat{r}')  
\Big]
\end{split}
\end{equation}
Likewise, a similar procedure leads to the representation of the dyadic Green's function of the magnetic field.
\begin{equation}
\label{equ:dyadGF_trace_M_exp}
\begin{split}
& \tilde {\bf G}^{tt}_{m0}(\bar{r},\bar{r}',t)=
 \sum_{n,m} \Big[  
 \dfrac{\partial_{r'}}{r'}  \big[ r' K_n^{(0)}(r,r',t) \big]   {\bf \Psi}_{nm}(\hat{r}) {\bf \Psi}_{nm}^{\ast}(\hat{r}') \\
& - \dfrac{\partial_{r}}{r} \big[ r K_n^{(0)}(r,r',t) \big]    {\bf \Phi}_{nm}(\hat{r}) {\bf \Phi}_{nm}^{\ast}(\hat{r}')
\Big]
\end{split}
\end{equation}
Whereas the above expressions are different from  those that are used conventionally, their scalar analogues have proven to be both robust (in terms of stability) and accurate.

\subsection{Volterra Integral kernels and Reduced Time-domain Integral Equations}\label{sec:kernels}

The expressions for the Green's functions (or its tangential trace) depend only on VSHs, and thanks to \eqref{eq:dyadGF_trace_def} they can be used instead of the dyadic Green's functions in the expressions for the EFIE and the MFIE. Specifically, it can be shown that by replacing the dyadic Green's function with its tangential trace, one can get the following integral operators; 
\begin{equation}
\mathbf{\mathcal{L}}^{tt}(\mathbf{X}(\bar{r},t)) =\mu \partial_t \int_{S'} \tilde{ \bf G}_{e0}^{tt}(\bar{r},\bar{r}') \otimes \mathbf{X}(\bar{r},t) dS'
\end{equation}
\begin{equation}
\mathbf{\mathcal{K}}^{tt}(\mathbf{X}(\bar{r},t)) = - \int_{S'}  \tilde {\bf G}_{m0}^{tt}(\bar{r},\bar{r}') \otimes \mathbf{X}(\bar{r},t) dS' 
\end{equation}
The above two operators are effectively equivalent to those in \eqref{eq:TDEFIE} and \eqref{eq:TDMFIE}, but permit the choice of vector spherical harmonics as basis functions for the unknown current densities.  Using the above equations together with Galerkin testing (that exploits the orthogonality of the VSHs) results in a set of one-dimensional Volterra equations. For the TD-EFIE, we obtain the following. 
\begin{subequations}
\begin{equation}
 \begin{split}
&  <{\bf \Psi}_{nm}^{\ast}({\hat r}), \mathbf{\mathcal{L}}( {\bf \Psi}_{n'm'}({\hat r}') )>  
=  <{\bf \Psi}_{nm}^{\ast}({\hat r}), \mathbf{\mathcal{L}}^{tt}( {\bf \Psi}_{n'm'}({\hat r}') )  > \\
& ~~~~~ = \dfrac{\mu c^2\partial_t^{-1}}{rr'}\partial_r\partial_{r'} \Big [  rr'K_n^{(0)}(r,r',t) \Big]
\delta_{nn'}\delta_{mm'} \\
 &~~~~~ \equiv K_n^{(1)}(r,r',t) \delta_{nn'}\delta_{mm'}
\end{split}
\end{equation}
\begin{equation}
 \begin{split}
& <{\bf \Phi}_{nm}^{\ast}({\hat r}), \mathbf{\mathcal{L}}( {\bf \Phi}_{n'm'}({\hat r}')  )  > =
<{\bf \Phi}_{nm}^{\ast}({\hat r}), \mathbf{\mathcal{L}}^{tt}( {\bf \Phi}_{n'm'}({\hat r}')  )  > \\
& ~~~~~~~ = \dfrac{\mu \partial _t}{c} \Big [ K_n^{(0)}(r,r',t)  \Big ]
\delta_{nn'}\delta_{mm'}   \\
& ~~~~~~~ \equiv K_n^{(2)}(r,r',t) \delta_{nn'}\delta_{mm'}
\end{split}
\end{equation}
\end{subequations}
where the functions $K_n^{(1)}$ and $K_n^{(2)}$ are the one-dimensional Volterra integral kernels for ${\bf \Psi}_{nm}$ and ${\bf \Phi}_{nm}$, respectively. The off-diagonal elements $\left<{\bf \Psi}_{nm}^{\ast}({\hat r}), \mathbf{\mathcal{L}}^{tt}( {\bf \Psi}_{n'm'}({\hat r}')  )\right> $ and $\left<{\bf \Phi}_{nm}^{\ast}({\hat r}), \mathbf{\mathcal{L}}^{tt}( {\bf \Phi}_{nm}({\hat r}') )\right>$ are identically zero. These equations can then be used to obtain the requisite Volterra integral equations for each of the two modes as 
\begin{subequations}
\begin{equation} \label{equ:VIE1}
 \begin{split}
K_n^{(1)}(r,r',t) \otimes J_{nm} ^1(t) & =<{\bf \Psi}_{nm}^{\ast}({\hat r}), \hat{n} \times \hat{n} \times {\bf E}^i (\bar{r},t)   >\\
& \equiv f_{nm}^1(t)
\end{split}
\end{equation}
\begin{equation} \label{equ:VIE2}
 \begin{split}
K_n^{(2)}(r,r',t) \otimes J_{nm}^2(t) & =<{\bf \Phi}_{nm}^{\ast}({\hat r}), \hat{n} \times \hat{n} \times {\bf E}^i (\bar{r},t)   >\\
& \equiv f_{nm}^2(t)
\end{split}
\end{equation}
\end{subequations}
In the above equations, the right-hand-sides are the projection of  $\hat{n} \times \hat{n} \times {\bf E}^i (\bar{r},t)$ onto the two VSH basis sets.

A similar procedure can be followed to derive the necessary equations for the TD-MFIE. The Volterra integral kernels $K_n^{(3)}$ and $K_n^{(4)}$ can be defined as follows.
\begin{subequations}
\begin{equation}
 \begin{split}
&  <{\bf \Psi}_{nm}^{\ast}({\hat r}), \mathbf{\mathcal{K}}( {\bf \Psi}_{n'm'}({\hat r}') )>  
=  <{\bf \Psi}_{nm}^{\ast}({\hat r}), \mathbf{\mathcal{K}}^{tt}( {\bf \Psi}_{n'm'}({\hat r}') )  > \\
& ~~~~~ =  -\dfrac{\partial_{r'}}{r'} \big[ r' K_n^{(0)}(r,r',t) \big]  \\
 &~~~~~ \equiv K_n^{(3)}(r,r',t) \delta_{nn'}\delta_{mm'}
\end{split}
\end{equation}

\begin{equation}
 \begin{split}
& <{\bf \Phi}_{nm}^{\ast}({\hat r}), \mathbf{\mathcal{K}}( {\bf \Phi}_{n'm'}({\hat r}')  )  > =
<{\bf \Phi}_{nm}^{\ast}({\hat r}), \mathbf{\mathcal{K}}^{tt}( {\bf \Phi}_{n'm'}({\hat r}')  )  > \\
& ~~~~~~~ = \dfrac{\partial_{r}}{r} \big[ r K_n^{(0)}(r,r',t) \big]   \\
& ~~~~~~~ \equiv K_n^{(4)}(r,r',t) \delta_{nn'}\delta_{mm'}
\end{split}
\end{equation}
\end{subequations}

Using these kernels, the resulting Volterra equations for the TD-MFIE are 
\begin{subequations}
\begin{equation} \label{equ:VIE3}
 \begin{split}
& J_{nm} ^1(t) + K_n^{(3)}(r,r',t) \otimes J_{nm}^1(t) \\
 & ~~~~~~~~~~~~~=<{\bf \Psi}_{nm}^{\ast}({\hat r}), \hat{n} \times {\bf H}^i (\bar{r},t)   >
 \equiv f_{nm}^3(t)
\end{split}
\end{equation}
\begin{equation} \label{equ:VIE4}
 \begin{split}
& J_{nm} ^2(t) + K_n^{(4)}(r,r',t) \otimes J_{nm}^2(t) \\
&~~~~~~~~~~~~~~ =<{\bf \Phi}_{nm}^{\ast}({\hat r}), \hat{n} \times {\bf H}^i (\bar{r},t)   >
\equiv  f_{nm}^4(t)
\end{split}
\end{equation}
\end{subequations}

Equations \eqref{equ:VIE1}, \eqref{equ:VIE2}, \eqref{equ:VIE3} and \eqref{equ:VIE4} constitute four equations that yield independent coefficients for the current on the surface of the sphere. In some sense, they can be considered a VSH transform of the integral operator. If necessary, one can combine the requisite integrals to form an effective time domain combined field solution as well. That said, it should be noted that the solution to scattering from a sphere has been reduced to solving these two equations that are uncoupled. In what follows, we shall provide a method for solving these equations numerically. 

\section{Time-dependent Debye-Mie series solution}\label{sec:time_mie}

Analytic solution to the above Volterra equations do not exist. As a result, numerical solution to these equations rely on discrete representation of $J_{nm}^1(t)$ and $J_{nm}^2(t)$. Note that, in this section, subscript ($nm$) for currents, coefficient vectors, right-hand-side vectors and system matrix representations is suppressed.  The approach that we follow is to partition the time axis uniformly, and represent basis and test function whose support reside on each time step. The basis functions can be any set of polynomials \cite{Pray2014}, however, we choose Legendre polynomials with support over the time step. Specifically, our representation 
\begin{equation}
\label{equ:td_exp1}
J^{\xi}(t_i + \tau\Delta_t) = \sum_{j=0}^{N_p} J^\xi_{ij} \Pi_i P_j(2\tau-1);~~~~~\tau\in[0,1)
\end{equation}
where $\Pi_i = 1$ for $t \in [t_i,t_{i+1})$, $P_j (\cdot)$ is a Legendre polynomial, $\xi = 1,2$, and $i = 1 \cdots N_t$ where $N_t$ is the total number of time steps. Using Galerkin testing results in the marching-on-in-time system. In stating the requisite equations, we use a block matrix representation that is indexed by time step as opposed to unknowns, viz., 
\begin{equation} \label{equ:mot}
{\cal Z}^{\nu}_{0} {\cal I}^\xi_{q} =   {\cal V}^{\nu}_{q} - \sum_{j=j_0}^{q-1} {\cal Z}^\nu_{q-j} {\cal I}^\xi_{j}
; ~~~~~\nu = 1, 2, 3, \text{or }  4
\end{equation}
where $j_0 = max(0,q-N_k)$, $q$ is the time step index ranging from $0$ to total number of steps $N_t$, and $N_k$ is the length of the discrete kernel. Here, ${\cal I}^\xi_{j}= \left [J^\xi_{j1}, \cdots, J^\xi_{jN_p}\right]^T $ are current coefficients for those polynomials associated with time step $j$. In a similar fashion, ${\cal Z}^{\nu}_{j}$ and ${\cal V}^{\nu}_{j}$ (the one-index subscript denotes the discritization in time) are ,respectively, discrete system and right-hand-side term at $j$th time step. At any given time $t$, the right-hand-side vector can be written as
\begin{equation}
\begin{split}
{\cal V}^\nu_{q} = \left [ \left <\Pi_q P_1 (\cdot), f^\nu (\cdot) \right >, \cdots, \left < \Pi_q P_{N_p} (\cdot) , f^\nu(\cdot) \right > \right ]^T
\end{split}
\end{equation}
and the element of the system matrix can be written as
\begin{equation} \label{equ:zmatrix}
{\cal Z}^{\nu}_{pi,qj} = < \Pi_p P_i (\cdot), K_n^{(\nu)}(\cdot)\otimes \Pi_q P_j(\cdot) > 
\end{equation}
In the expressions for vector and matrix elements, the two-index and four-index subscripts denote the discritization in space.

In the above equations, the temporal dependencies of the function follows from their explicit definitions earlier in the Section. It should be noted that $K_n^{(1)}$ has an infinite tail, very similar to that encountered while evaluating the contribution of the scalar potential in regular TD-EFIE \cite{Pray2014}. As in \cite{Pray2012,Shanker2003}, introduction of an auxiliary charge ameliorates the computational complexity from ${\cal O}(N_t^2) $ to $ {\cal O}(N_t)$. Introducing an auxiliary charge the MOT system for $K_n^{(1)}$ can be written as 
\begin{equation} \label{equ:mot_charge}
{\cal Z}^{I}_{0} {\cal I}^1_{q} =   {\cal V}^{1}_{q} -\sum_{j=j_0}^{q-1} {\cal Z}^I_{q-j} {\cal I}^1_{j} 
-  {\cal C}_{j-1}
\end{equation}
where $j_0 = max(0,q-N_k)$, $N_k$ is the length of the non-constant part of the discrete kernel $K_n^{(1)}$. The auxiliary charge coefficient is 
chosen as
\begin{equation}
\begin{split}
 & C_j= C_{j-1} +  \sum_{k=k_0}^{j}{\cal Z}^I_{j-k} {\cal I}_k
\end{split}
\end{equation}
where $k_0 = \min(0, j-N_k)$.
The elements of the modified system (${\cal Z}^{I}$)
are defined as follows.
 \begin{equation} \label{equ:zmatrix_I}
{\cal Z}^{I}_{k} = {\cal Z}^{1}_{k} -{\cal Z}^{1}_{k-1} 
\end{equation}
Equations \eqref{equ:mot} and \eqref{equ:mot_charge} can be solved for the coefficients in every time step, and of course for each harmonic. On another brief note regarding discretization of these equation; since both testing and source basis functions are pulse functions that are multiplied by a higher order polynomial, it can be shown that the inner product in \eqref{equ:zmatrix} can be equivalently written as a convolution of the kernel $K_n^{(\nu)}$ with an interpolatory polynomial that is higher order with support from $(-\Delta_t,\Delta_t)$ and tested by a delta function. As a result, they are satisfy the function spaces prescribed in \cite{TERRASSE1993}. The resulting matrix elements are smooth, well behaved, and causal. More importantly, they do not require anterpolation unlike that used in \cite{Li2014,Weile2004}. In what follows, we briefly describe the eigen-spectrum 
analysis \cite{Walker2002,Pray2012} of the kernels presented for sample harmonics. 

\section{Eigen-spectrum analysis of the MOT system}

Next, we examine the stability of the MOT system developed earlier by analyzing the eigen-spectrum of the marching system for different harmonics. We note that this section is provided purely for completeness; the equations provided below are similar to those presented in \cite{Walker2002,Pray2012} by comparing the similarity of the kernels analyzed herein to those presented in these papers.  Specifically, for equations involving $K_n^{(1)}$, the MOT can be  rewritten in the following matrix system form
\begin{eqnarray} \left[ \begin{array}{c|c}
  \underline{A}_{11} & \underline{0} \\
  \underline{0} & \underline{\mathbb I} \\
  \end{array}\right] \underline{{ I}}_{j+1} = \underline{{F}}_{j+1}
  -  \left[ \begin{array}{c|c}
  \underline{B}_{11} & \underline{B}_{12} \\
  \underline{B}_{21} &\underline{B_{21}}  \\
  \end{array}\right]\underline{{I}}_{j}, 
\label{eq:marchmat_z}
  \end{eqnarray}
where $\underline{{ I}}_j = [I_j^T,\cdots, J_{j-N_k}^T, C_{j-1}^T,\cdots, C_{j-1-N_k}^T ]^T$, 
$\underline{{ F}}_j = [I_j^T,0, \cdots, 0, 0,\cdots, 0 ]^T$ (here the superscript $T$ denotes transpose)
and $\underline{\mathbb I} $ is the indentity matrix.
The elements in other matrices can be computed as follows.
\begin{subequations}
\begin{eqnarray}
  \underline{A}_{11}=&\left[ \begin{array}{ccccc}
\left[{\mathcal Z}^I_0\right] & 0 &   & &   \ldots \\
0 &\left[{\mathbb I}\right] & 0 &    & \ldots \\
0 & 0 &\left[{\mathbb I}\right] & 0   & \ldots \\
. & & &  & \\
\ldots & & &  & \left[{\mathbb I}\right]
  \end{array}\right]~,\\
  \underline{B}_{11}=&\left[ \begin{array}{cccc}
\left[{\mathcal Z}^I_1\right] & \left[{\mathcal Z}^I_2\right] & \ldots &  \left[{\mathcal Z}^I_{N_k}\right] \\
  \left[{\mathbb I}\right] &0 &    & \ldots \\
   0 &\left[{\mathbb I}\right]  &  & \ldots \\
 . & &  & \\ 
 . & &  & \\
 \ldots & &\left[{\mathbb I}\right]&0 
  \end{array}\right]~,\\
  \underline{B}_{12}=&\left[ \begin{array}{ccccc}
 \left[{\mathbb I}\right] &0 &  &    & \ldots  \\
   0 &0 &  &   & \ldots \\
   . & &  & & \\
   . & &  & & \\
   0 & 0 & & & \ldots
   \end{array}\right]~,\\
  \underline{B}_{21}=&\left[ \begin{array}{cccc}
\left[{\mathcal Z}^I_1\right] & \left[{\mathcal Z}^I_2\right] & \ldots &  \left[{\mathcal Z}^I_{N_k}\right]\\
 0 &0  &  & \ldots\\
    . & & & \\
   . & & & \\
    \ldots & &0 &0  \\
   \end{array}\right]~,\\
  \underline{B}_{22}=&\left[ \begin{array}{ccccc}
 \left[{\mathbb I}\right] & 0 &  &    & \ldots\\
  \left[{\mathbb I}\right] & 0 &  &    & \ldots\\
  . & &  & & \\
  . & &  & &\\
  \ldots &  & &\left[{\mathbb I}\right]&0
  \end{array}\right]~,
\end{eqnarray}
\end{subequations}

For systems involving other kernels, the corresponding set of equations are
\begin{equation}
 \underline{A} ~ \underline{{ I}}_j = \underline{{ F}}_j -  \underline{B} ~\underline{{ I}}_{j-1}
\end{equation}
where $\underline{{ I}}_j = [I_j^T,\cdots, J_{j-N_k}^T]^T$, 
$\underline{{ F}}_j = [I_j^T,0, \cdots, 0 ]^T$. Compared to the updating equations for the first kernel, only current vectors are considered, and the elements of the matrices $A$ and $B$ can be similarly obtained by changing the corresponding kernels in $ \underline{A}_{11}$ and $ \underline{B}_{11}$. 

To obtain stability of the marching system, eigen-spectrum analysis should be done for matrix $ \underline{A}^{-1} \underline{B}$. And in next section, the results of this analysis are presented for different kernels and harmonics.

\section{Numerical Examples}\label{sec:results}
The contribution of this work mainly lies in the derivation reduced time-domain integral equation for transient Debye-Mie series and its marching on in time solution. Hence we provide several examples of using TD-EFIE and TD-MFIE are given to test the accuracy of the method and to demonstrate the stability of the MOT algorithm. Because the spatial spectrum is the same for both time and frequency domain and the orthogonality allows the mode-by-mode solution. In what follows, comparison for a single mode is presented as opposed to their sum.
The problem studied here is scattering  by a perfectly electrically conducting (PEC) sphere of radius $a=1$m. 
The incident field is a plane wave with a modulated Gaussian profile
\begin{equation}
	{\bf E}^{i}({\bf r},t)={\hat x}\cos(2\pi f_0t)e^{-(t-{\bar r}\cdot{\hat k} /c-t_p)^2/2\sigma^2 }
	\label{einc}
\end{equation}
where ${\hat k} =\hat{z}$ is the direction of propagation, $f_0$ is center frequency. In the expression , $\sigma=3/(2\pi B)$ and $t_p=40\sigma$, where $B$ denotes the bandwidth.  In all the numerical examples of this section,  center frequency of $f_0 = 0.4$GHz and bandwidth of $B=0.3$GHz are chosen. With this configuration, the frequency band we are interested in is between $0.1$GHz and $0.7$GHz. Due to the fact that the power spectral density value at $f_{max}=0.7$GHz is $39$dB lower than that at the center frequency, the number of spherical harmonics to (spatially) represent the plane wave is bounded by $N_m\approx [2ka]  = 30 $.
\begin{figure} 
   \centering
   \includegraphics[width=1\columnwidth,keepaspectratio=true]{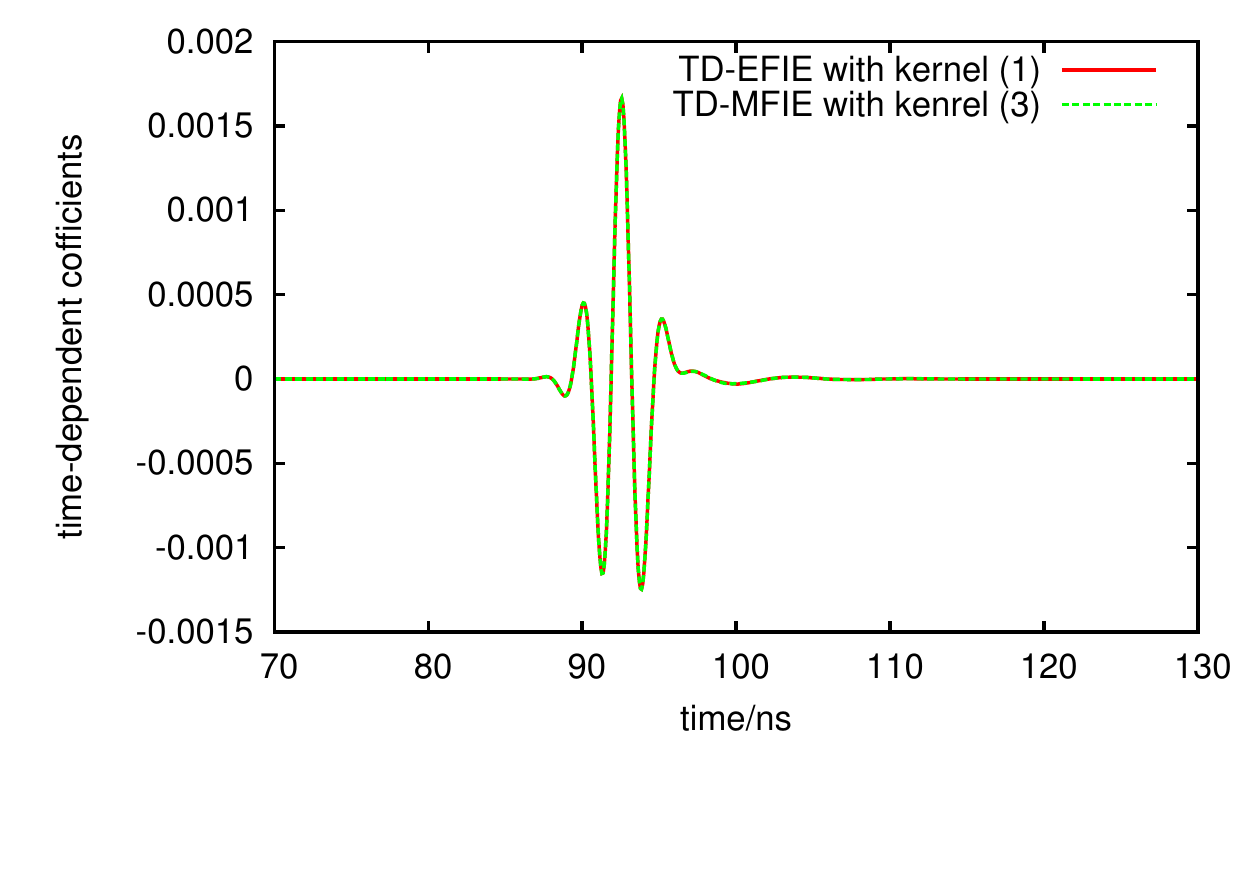}
\caption{Time-domain coefficients for mode ${\bf \Psi}_3^{1}$}
\label{fig:fig1}
\end{figure}

\begin{figure}
   \centering
   \includegraphics[width=1\columnwidth,keepaspectratio=true]{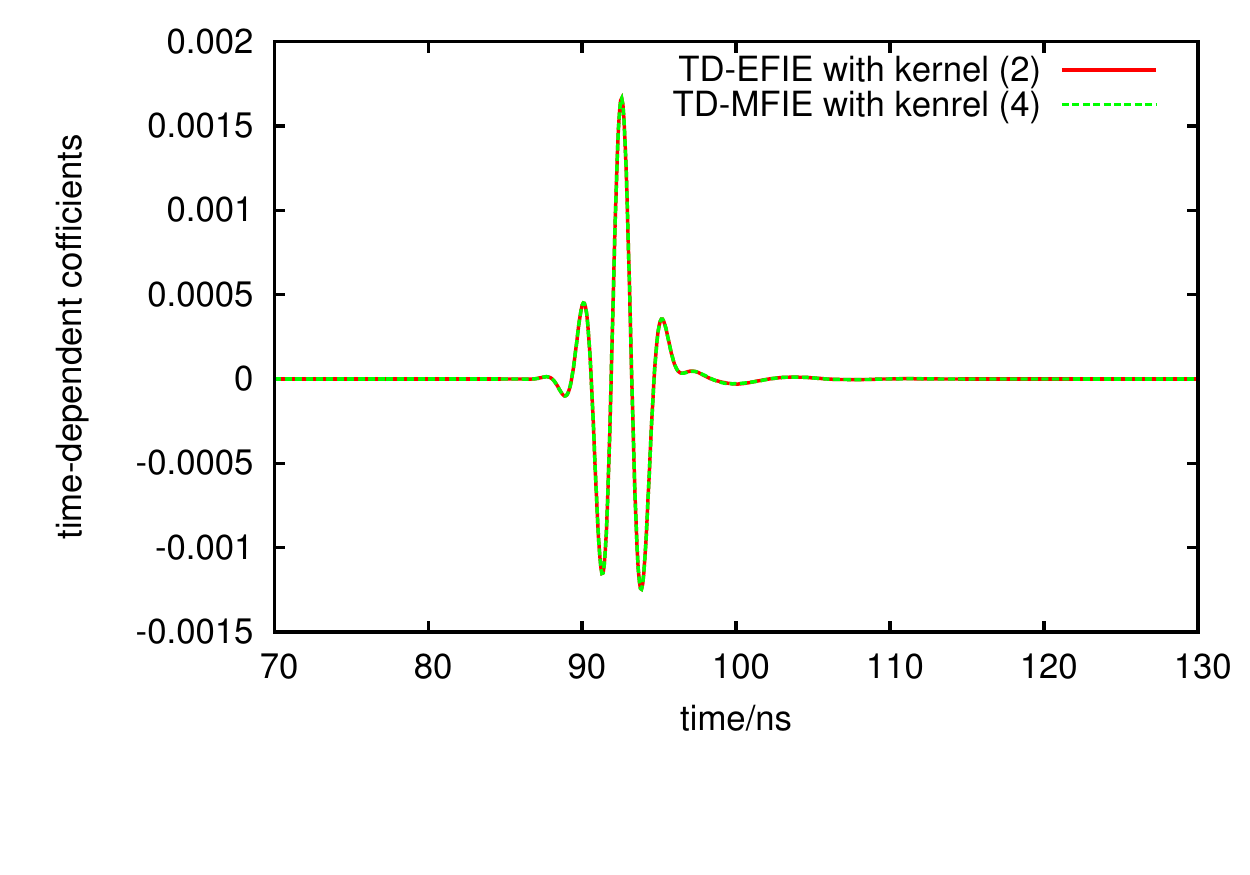}
\caption{Time-domain coefficients for mode ${\bf \Phi}_{3}^{1}$}
\label{fig:fig2}
\end{figure}

\begin{figure}
   \centering
   \includegraphics[width=1\columnwidth,keepaspectratio=true]{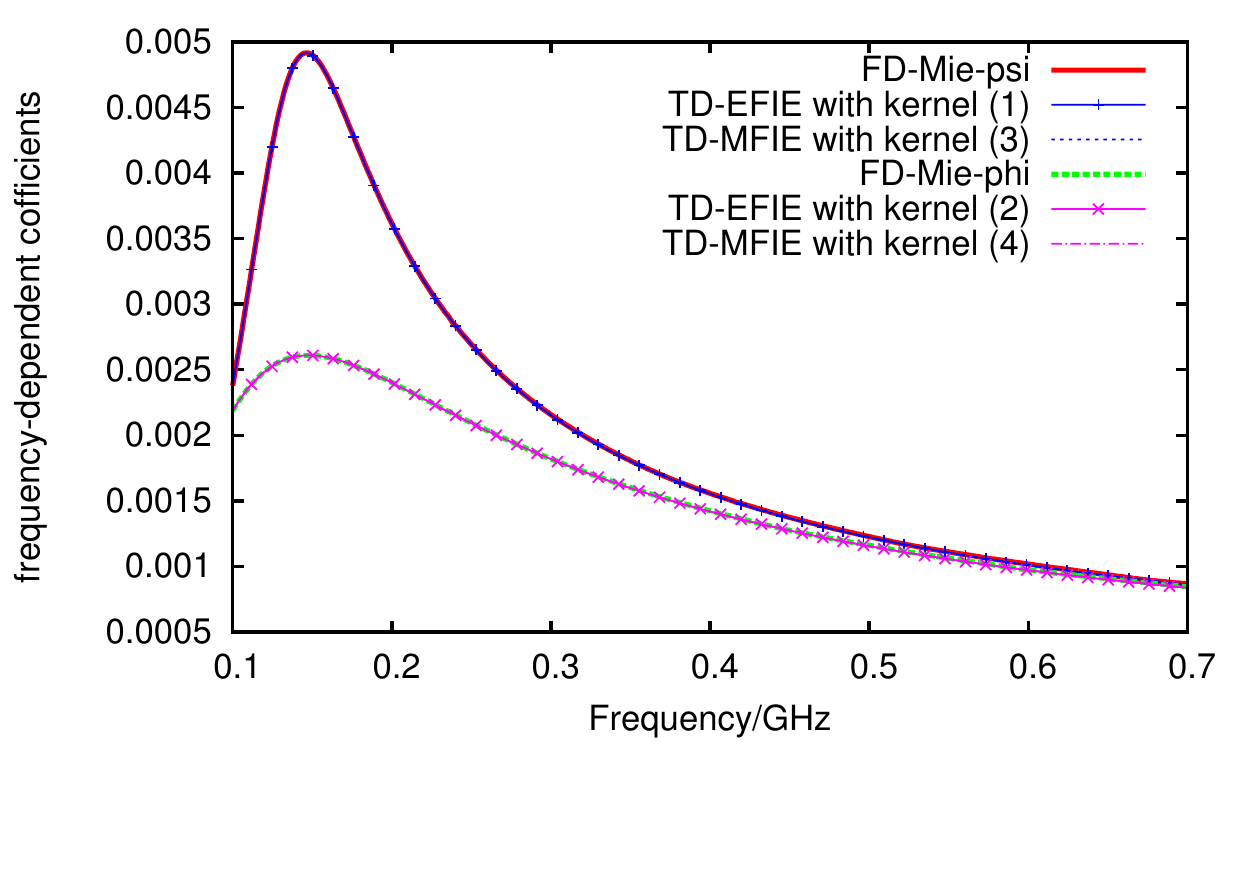}
\caption{Frequency-domain coefficients for modes ${\bf \Psi}_3^{1}$ and  ${\bf \Phi}_{3}^{1}$,compared to two curves (denoted by FD-Mie-psi and FD-Mie-phi respectively) with frequency domain Mie seires approach. }
\label{fig:fig3}
\end{figure}

\begin{figure} 
   \centering
   \includegraphics[width=1\columnwidth,keepaspectratio=true]{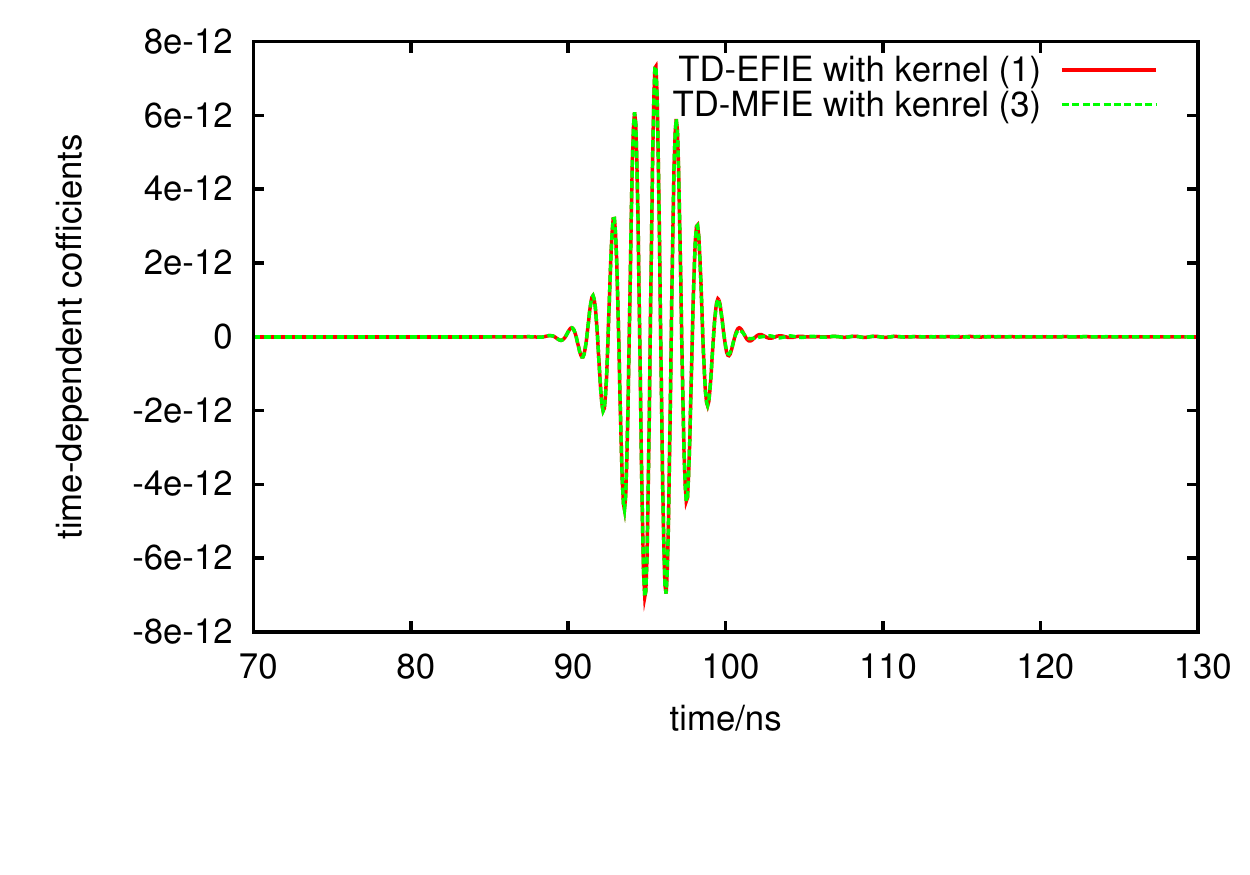}
\caption{Time-domain coefficients for mode ${\bf \Psi}_{30}^{1}$}
\label{fig:fig4}
\end{figure}

\begin{figure}
   \centering
   \includegraphics[width=1\columnwidth,keepaspectratio=true]{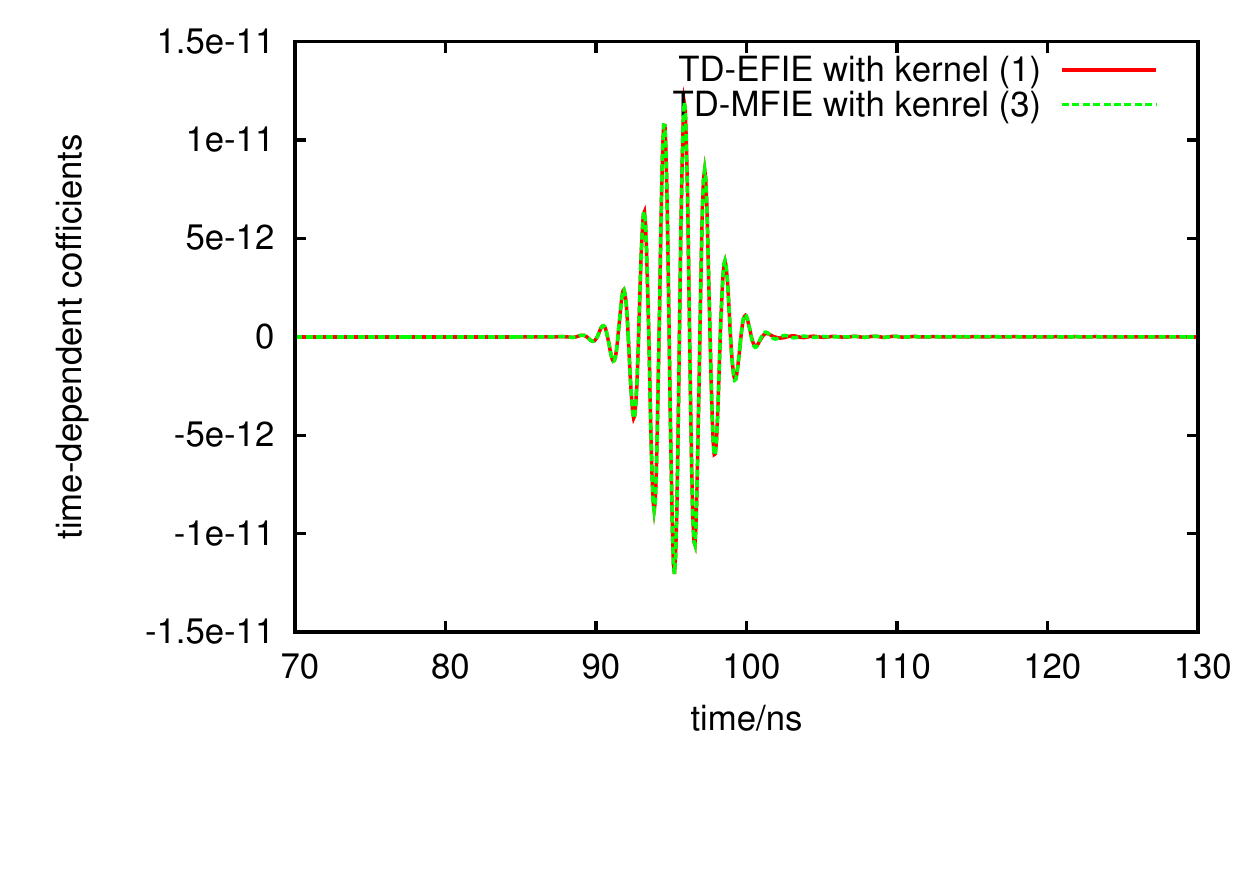}
\caption{Time-domain coefficients for mode ${\bf \Phi}_{30}^{1}$}
\label{fig:fig5}
\end{figure}

\begin{figure}
   \centering
   \includegraphics[width=1\columnwidth,keepaspectratio=true]{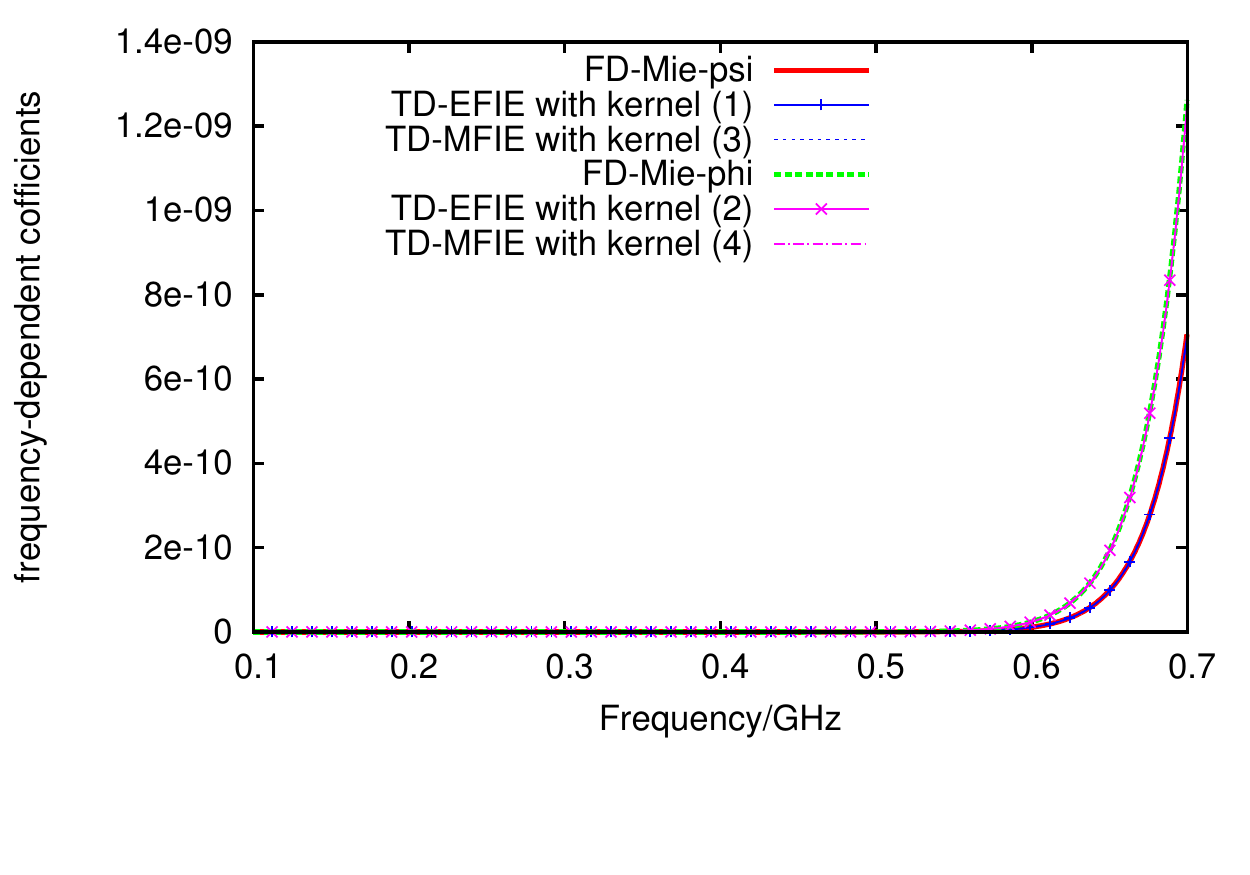}
\caption{Frequency-domain coefficients for mode ${\bf \Psi}_{30}^{1}$ and mode ${\bf \Phi}_{30}^{1}$,,compared to two curves (denoted by FD-Mie-psi and FD-Mie-phi respectively) with frequency domain Mie seires approach. }
\label{fig:fig6}
\end{figure}

The first test is to illustrate the accuracy of the proposed method by comparing its frequency-domain counterpart. In this test, the temporal basis function are chosen up to first order Legendre polynomials. The time-dependent coefficients for modes ${\bf \Psi}_3^{1}$ and ${\bf \Phi}_{3}^{1}$ are plotted in Fig. \ref{fig:fig1} and Fig. \ref{fig:fig2}. Each plot contains two curves that are obtained using both TD-EFIE and TD-MFIE with corresponding kernels. The comparison in frequency domain within the frequency band is demonstrated in Fig.\ref{fig:fig3}, and the results from frequency domain Mie (FD-Mie) series are also provided to verify the frequency responses. One can observe that the agreement is excellent between time and frequency domain. Similarly conclusion can be drawn for higher order modes. Time-domain and frequency-domain results for modes ${\bf \Phi}_{30}^{1}$ and ${\bf \Phi}_{30}^{1}$  are given in Figs. \ref{fig:fig4}-\ref{fig:fig5} and Fig. \ref{fig:fig6}, respectively. One can also observe 
that  coefficients for modes of degree 30 are orders of magnitude smaller that those of degree 3. From these results, one can 
recognize the feasibility of extracting transient response of spherical object with the presented mode by mode MOT approach. In the simulations we conducted, relative error convergence in frequency domain (down to $10^{-12}$) versus the order of the temporal basis function (up to 6) has been observed.

\begin{figure}
   \centering
   \includegraphics[width=1\columnwidth,keepaspectratio=true]{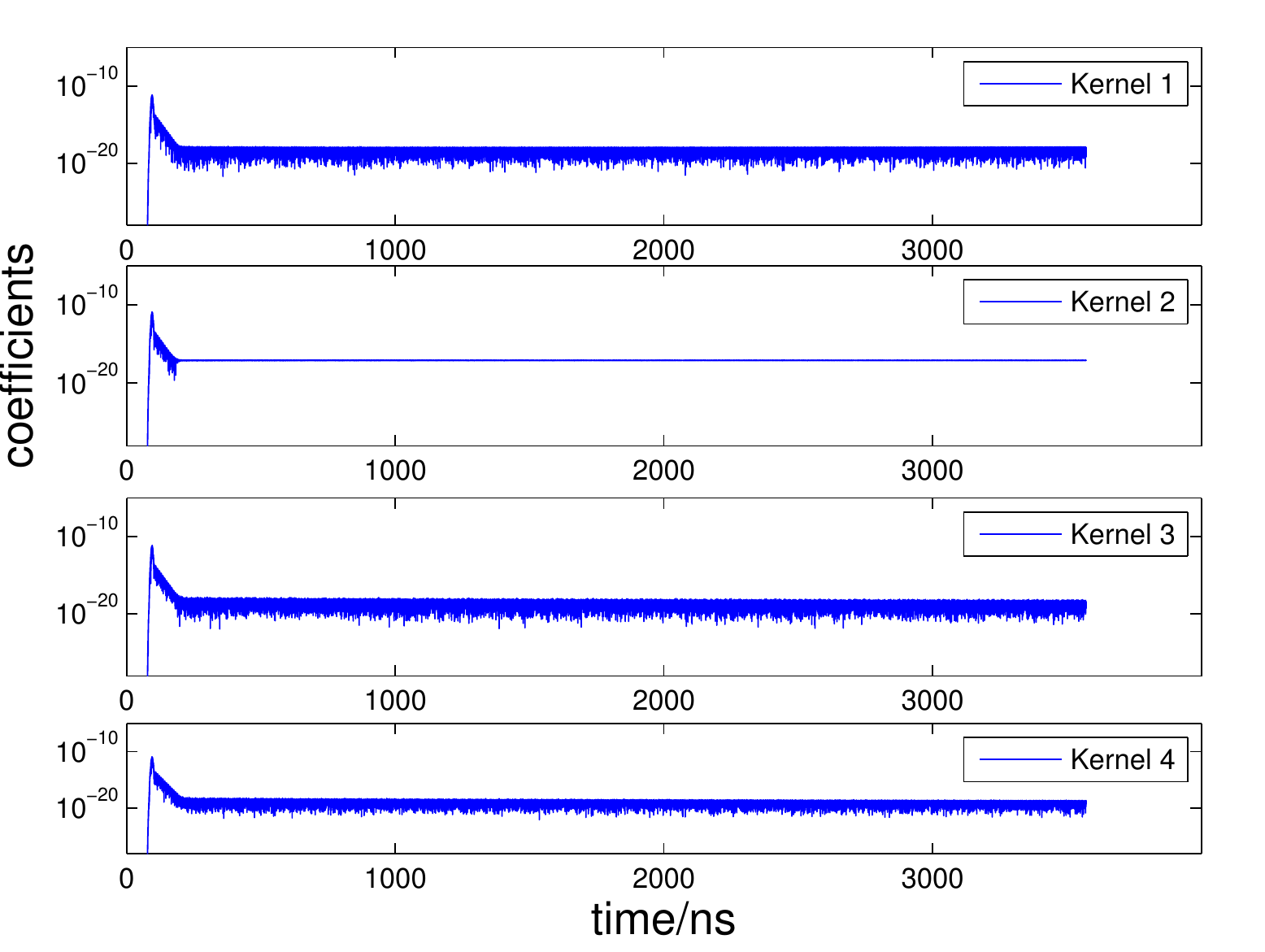}
\caption{Late-time stability}
\label{fig:stable}
\end{figure}

\begin{figure}
  \centering
\subfloat[]{
 \includegraphics[width=0.5\columnwidth,keepaspectratio=true]{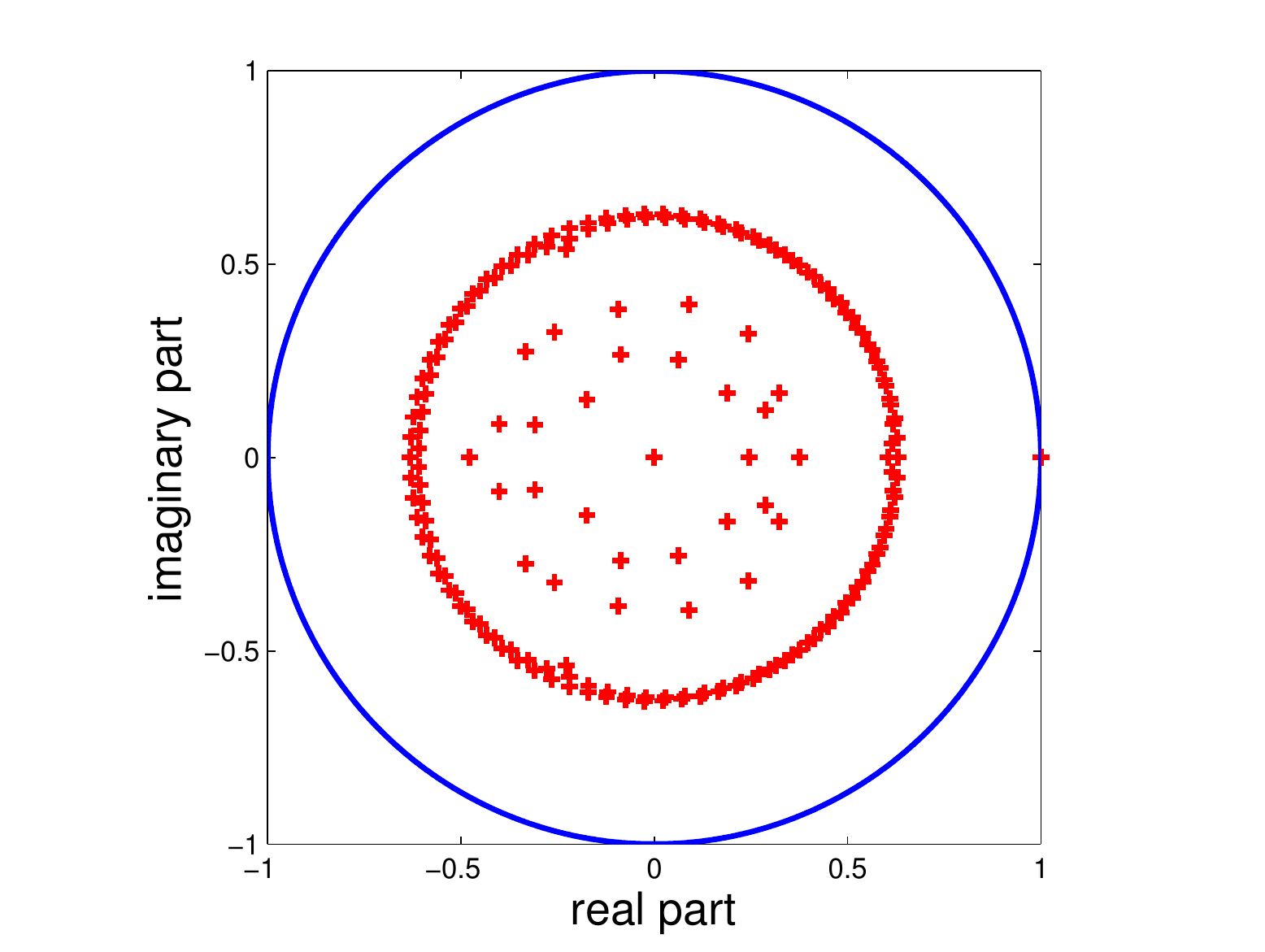}
                \label{fig:eigen_k1}
 }
\subfloat[]{
 \includegraphics[width=0.5\columnwidth,keepaspectratio=true]{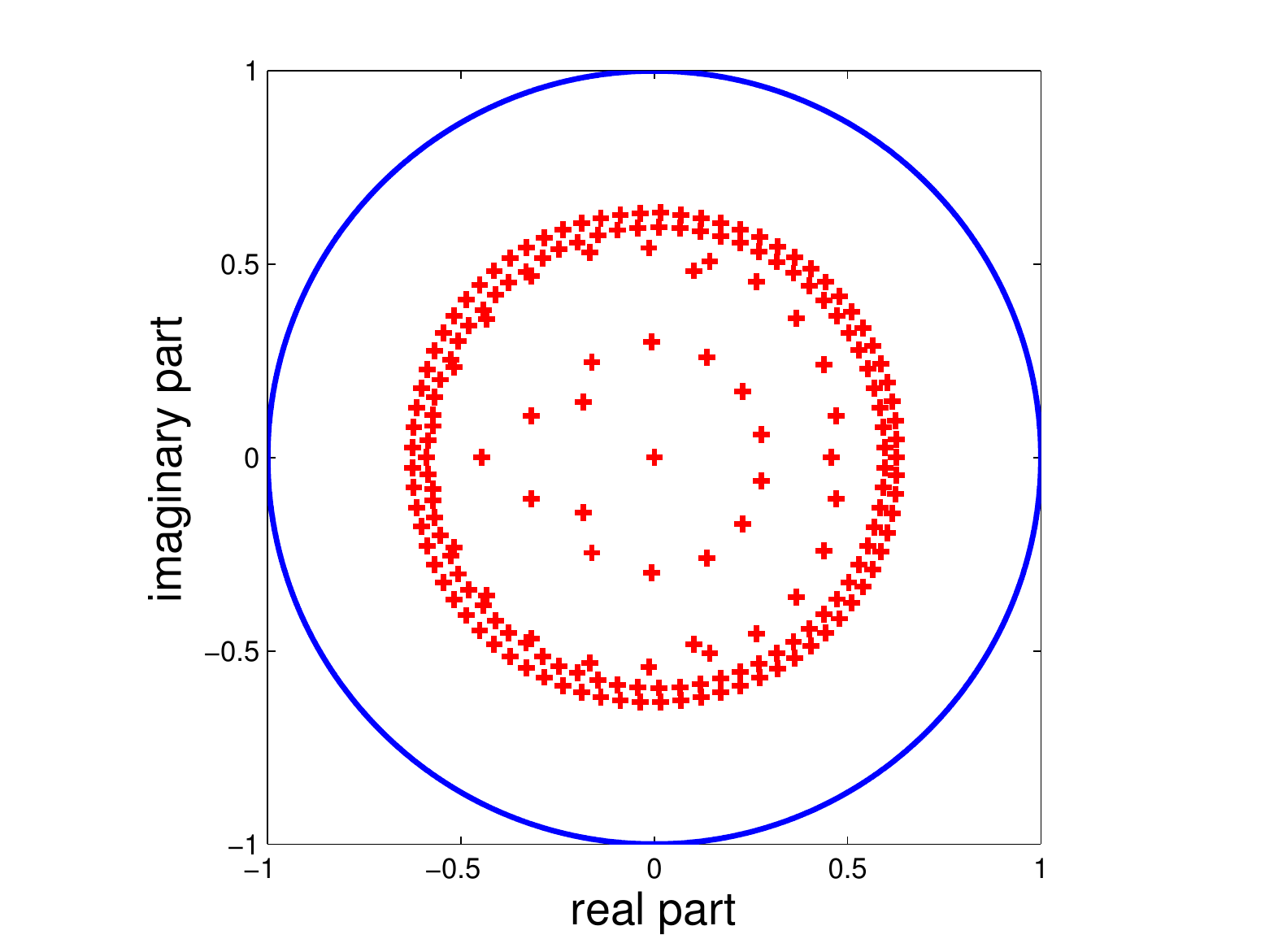}
                \label{fig:eigen_k2}
 }

\hfill

\subfloat[]{
 \includegraphics[width=0.5\columnwidth,keepaspectratio=true]{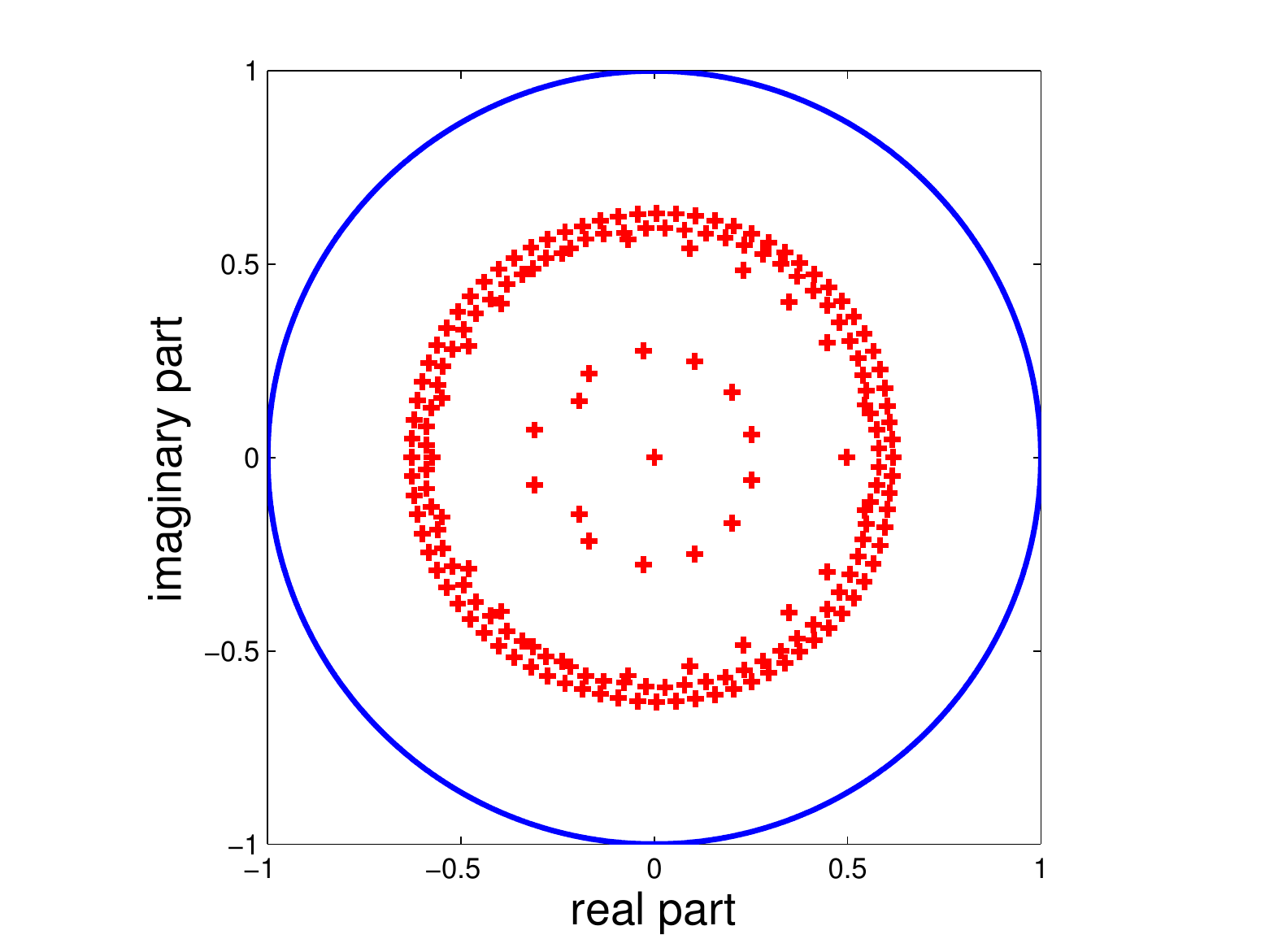}
                \label{fig:eigen_k3}
 }
\subfloat[]{
 \includegraphics[width=0.5\columnwidth,keepaspectratio=true]{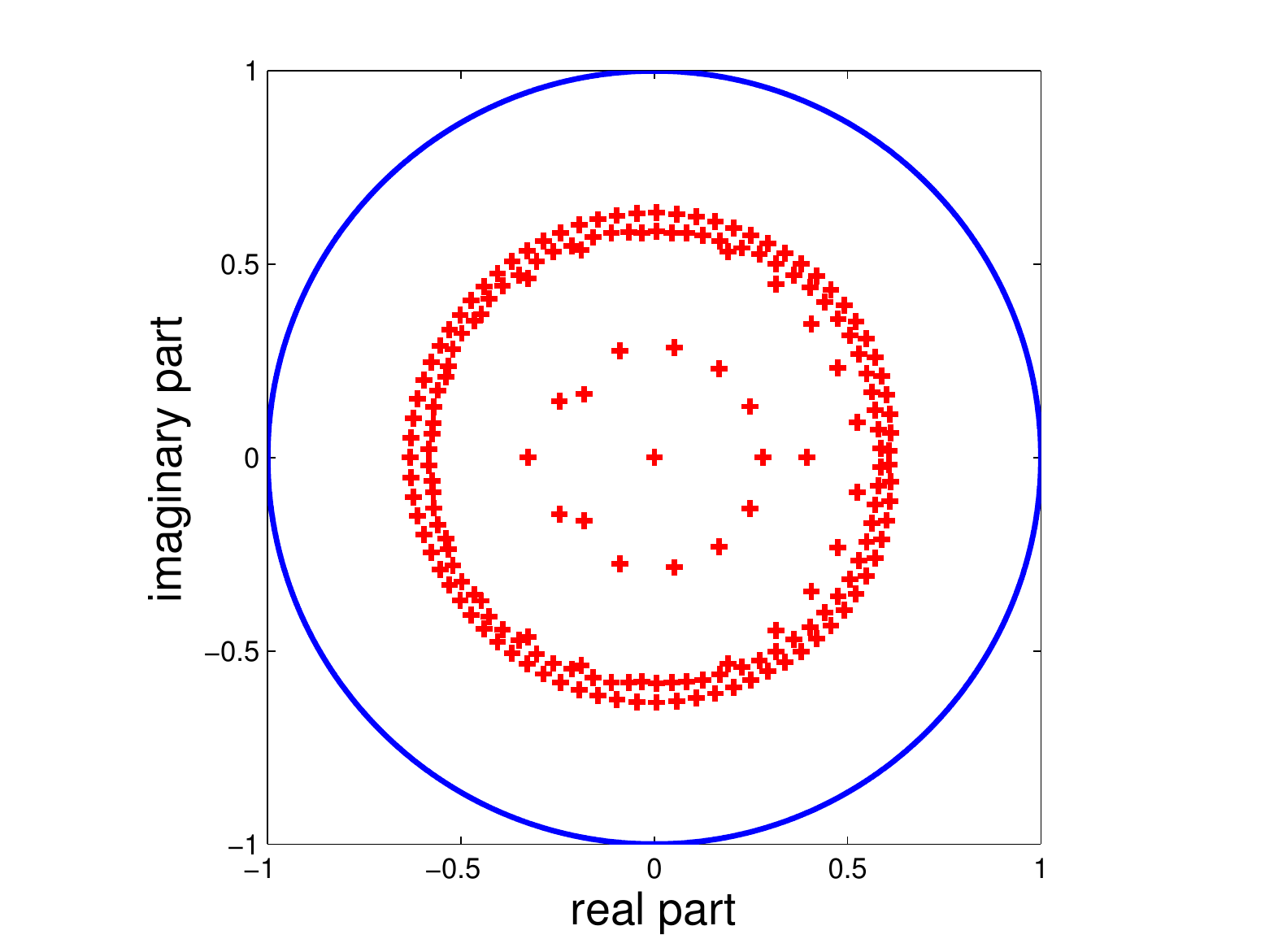}
                \label{fig:eigen_k4}
 }
   \caption{Eigenvalues of MOT system of (a) kernel 1, (b) kernel 2, (c) kernel 3 and (d) kernel 4.} 
\label{fig:eigen}
\end{figure}

The second test is to demonstrate the stability of the marching algorithm. Four 100,000-step simulations corresponding to the derived four reduced Volterra kernels are carried out, where the temporal basis functions are chosen up to first order. The solutions are transient coefficients for the modes  ${\bf \Psi}_{30}^{1}$ (solutions to VIE of kernels 1 and 3)  and ${\bf \Phi}_{30}^{1}$ (solutions to VIE of kernels 2 and 4). The time-domain coefficients are plotted in Fig. \ref{fig:stable}. It's obvious that no late-time instability is observed in these simulations. To further show the stability of the resulting MOT scheme, the eigenvalues analysis of the discretized marching system outlined earlier is given. The real and imaginary parts of the eigenvalues are plotted in Figs. \ref{fig:eigen_k1}-\ref{fig:eigen_k4}, with all the values are within the unit circle except one lying on the unit circle which is associated with kernel 1. Similar stability can be also found for simulations with higher order temporal 
basis functions and Volterra kernels of different degrees.


\section{Conclusions}\label{sec:conc}
In this work, a TDIE-based direct approach is proposed to calculate the time-dependent spherical multipoles due to the scattering by a PEC sphere. By using spherical harmonics as basis functions and expanding the tangential trace of time-domain dyadic Green's function with VSH, the proposed method analytically reduces the original time-domain integral into several novel one-dimensional Volterra integral equations. Those reductions allow very efficient computation for those time-domain spherical multipoles in a marching manner. As shown in the numerical examples, the accuracy and convergence of discretization in time are verified. The resulting kernels could also be considered as good testbeds for various temporal basis function choices and testing schemes. Stability issue arising in conventional TDIE could be studied with those reduced integral equation kernels, which is a current research topic. The presented technique could also be used for multiple particles and enhanced by fast methods\cite{
Shanker2003}, which is another part of our future work.

\section*{Acknowledgement}

The authors wish to acknowledge the HPCC facility at Michigan State University, support from NSF CCF-1018516 and NSF CMMI-1250261, and support from the Department of Defense CREATE-RF program. B. Shanker is grateful for support from the Department of Defense High Performance Computing Modernization Program Office (HPCMPO) and its Computational Research and Engineering for Acquisition Tools and Environments (CREATE) Program.

\section*{Appendix}
\subsection{Vector Spherical Harmonics (VSH) and Wave Functions (VSWF)}
The two vector spherical wave functions used in conventional expansion of dyadic Green's function are defined as follows.
\begin{equation}
  {\bf M}_{nm}^{(i)}(k,\bar{r}) =  - z_n^{(i)}(kr)  \mathbf{\Phi}_{nm}(\theta,\phi) 
\end{equation}

\begin{equation}
\begin{split}
  {\bf N}_{nm}^{(i)}(k,\bar{r}) =  & \dfrac{[kr z_n^{(i)}(kr) ]'}{kr} \mathbf{\Psi}_{nm}(\theta,\phi) \\
  &  + \dfrac{\sqrt{n(n+1)}}{kr}  z_n^{(i)}(kr) {\bf Y}_n^m(\theta,\phi) 
\end{split}
\end{equation}
where $ {\bf Y}_n^m(\theta,\phi) =  Y_n^m(\theta,\phi)\hat{r}$.

In order to derive the spherical harmonics expansion of the tangential trace of the dyadic Green's function, the mapping relations between VSH and VSWF should be used.
\begin{equation}
    \hat{r} \times {\bf M}_{nm}^{(i)}(k,\bar{r}) 
= z_n^{(i)}(kr)  \mathbf{\Psi}_{nm}(\theta,\phi)
\end{equation}

\begin{equation}
  \hat{r}  \times   \hat{r} \times {\bf M}_{nm}^{(i)}(k,\bar{r}) 
=  z_n^{(i)}(kr)  \mathbf{\Phi}_{nm}(\theta,\phi)
\end{equation}

\begin{equation}
\begin{split}
     \hat{r} \times {\bf N}_{nm}^{(i)}(k,\bar{r})   
 = \dfrac{[kr z_n^{(i)}(kr) ]'}{kr}   \mathbf{\Phi}_{nm}(\theta,\phi) 
\end{split}
\end{equation}

\begin{equation}
\begin{split}
     \hat{r} \times \hat{r} \times  {\bf N}_{nm}^{(i)}(k,\bar{r})    = -\dfrac{[kr z_n^{(i)}(kr) ]'}{kr}   \mathbf{\Psi}_{nm}(\theta,\phi) 
\end{split}
\end{equation}
\subsection{Spherical expansion of the time domain magnetic dyadic Green's function}
Frequency-domain dyadic Green's function of magnetic field is the solution to the wave equation.
\begin{equation}
 \nabla \times \nabla \times \tilde {\bf G}_{m0} - k^2 \tilde {\bf G}_{m0} =  \nabla \times [{\it \tilde I} \delta (\bar{r}-\bar{r}')],
\end{equation}

The dyadic Green's function can be written in terms of the scalar Green's function in free space.
\begin{equation}
\tilde {\mathbf G}_{m0}(\bar{r},\bar{r}',\omega)
= \nabla \times [ {\it \tilde I} G_o(\bar{r},\bar{r}') ]
=\nabla G_0(\bar{r},\bar{r}') \times {\it \tilde I}
\end{equation}

As in electric field case, the dyadic Green's function can be expanded using a series of vector spherical wave functions.

\begin{equation}
\begin{split}
 \tilde {\bf G}_{m0}(\bar{r},\bar{r}',\omega)
& = jk^2     \sum_{n,m} \Big[  
  {\bf M}_{nm}^{(4)}(k,\bar{r}) {\bf N}_{nm}^{(1)\ast}(k,\bar{r}') \\
 &~~~~~~~~~~~ +   {\bf N}_{nm}^{(4)}(k,\bar{r}) {\bf M}_{nm}^{(1)\ast}(k,\bar{r}')
\Big]
\end{split}
\end{equation}

By taking the tangential trace of this dyadic Green's function, one can get a modified Green's function expanded purely by vector spherical harmonics, separating the explicit radial dependence from non-radial dependence.
\begin{equation}
\begin{split}
&  \tilde {\bf G}_{m0}^{tt}(\bar{r},\bar{r}',\omega)
 = {\hat r} \times \tilde {\bf G}_{m0}(\bar{r},\bar{r}',\omega)  \times {\hat r}'  \times {\hat r}' \\
 & = jk^2 \sum_{n,m} \Big[  
  \dfrac{\big[  kr{z_n^{(4)}}(kr) \big]'}{kr} z_n^{(1)\ast}(kr') {\bf \Phi}_{nm}(\hat{r}) {\bf \Phi}_{nm}(\hat{r}') \\
 & ~~~~~~  - z_n^{(4)}(kr) \dfrac{\big[ kr'{z_n^{(1)\ast}}(kr') \big]'}{kr'}   {\bf \Psi}_{nm}(\hat{r}) {\bf \Psi}_{nm}(\hat{r}')
\Big]
\end{split}
\end{equation}
The reason that ${\hat r} \times$ instead of ${\hat r} \times{\hat r} \times$ is chosen in the left side is due to the fact that there is already one tangential operator ${\hat r} \times$ in the ${\bf { \cal K}}$ operator.
Similar as in electric field case, the modified dyadic Green's function has the following properties.
\begin{subequations}
 \begin{equation}
\tilde {\bf G}_{m0}^{tt}(\bar{r},\bar{r}',\omega) \cdot {\bf X}
=  - \tilde {\bf G}_{m0} (\bar{r},\bar{r}',\omega) \cdot {\bf X} 
\end{equation}
 \begin{equation}
 {\bf X} \cdot  {\hat r} \times  \tilde {\bf G}_{m0}^{tt}(\bar{r},\bar{r}',\omega) 
= - {\bf X} \cdot  {\hat r} \times   \tilde {\bf G}_{m0} (\bar{r},\bar{r}',\omega) 
\end{equation}
\end{subequations}

In order to get the time-domain expansion, inverse Fourier transform of the above equation gives the following formula.
\begin{equation}
\label{equ:dyadGF_M_exp}
\begin{split}
&  \tilde {\bf G}_{m0}^{tt}(\bar{r},\bar{r}',t)
  = \\
& \sum_{n,m} \Big[  
  {\cal F}^{-1}  \big\{  jk \dfrac{\partial}{r \partial r}\big[ r{z_n^{(4)}}(kr) \big] z_n^{(1)}(kr')   \big\} {\bf \Phi}_{nm}(\hat{r}) {\bf \Phi}_{nm}(\hat{r}') \\
 &   - {\cal F}^{-1}  \big\{  jk \dfrac{\partial}{r' \partial r'}\big[ r'{z_n^{(1)}}(kr') \big] z_n^{(4)}(kr)   \big\} 
  {\bf \Psi}_{nm}(\hat{r}) {\bf \Psi}_{nm}(\hat{r}')
\Big]
\end{split}
\end{equation}

Due to the linearity of (inverse) Fourier transform, one can get the following results.
\begin{equation}
\begin{split}
  {\cal F}^{-1}  \big\{  jk \dfrac{\partial}{r' \partial r'}\big[ r'{z_n^{(1)}}(kr') \big] z_n^{(4)}(kr)   \big\}  
 =  - \dfrac{\partial_{r'}}{r'} \big[ r' K_n^{(0)}(r,r',t) \big] 
\end{split}
\end{equation}
\begin{equation}
\begin{split}\\
  {\cal F}^{-1}  \big\{  jk \dfrac{\partial}{r \partial r}\big[ r{z_n^{(4)}}(kr) \big] z_n^{(1)}(kr')   \big\}
 = - \dfrac{\partial_{r}}{r} \big[ r K_n^{(0)}(r,r',t) \big] 
\end{split}
\end{equation}


By using above two kernels in \eqref{equ:dyadGF_M_exp}, one can get the spherical expansion for the time-domain dyadic Green's function as in \eqref{equ:dyadGF_trace_M_exp}.

%

\bibliographystyle{ieeetr}
\bibliography{tdexp}

\end{document}